\documentclass[12pt]{article}

\font\tenmsb=msbm10
\font\sevenmsb=msbm7
\font\fivemsb=msbm5

\newfam\msbfam
\textfont\msbfam=\tenmsb
\scriptfont\msbfam=\sevenmsb
\scriptscriptfont\msbfam=\fivemsb
\def\Bbb#1{{\fam\msbfam #1}}

\makeatletter
\ifnum\@ptsize=0  \fi \ifnum\@ptsize=2
 \fi 

\newcommand\qed{{\hspace*{\fill}Q.E.D.\vskip12pt plus 1pt}}

\newcommand\sE{{\cal E}}
\newcommand\sF{{\cal F}}
\newcommand\sG{{\cal G}}

\newcommand\sI{{\cal I}}

\newcommand\sO{{\cal O}}
\newcommand\sS{{\cal S}}

\newcommand\sC{{\cal C}}

\newcommand\sT{{\cal T}}

\newcommand\bR{{\Bbb R}}
\newcommand\bZ{{\Bbb Z}}
\newcommand\bC{{\Bbb C}}
\newcommand\bQ{{\Bbb Q}}
\newcommand\bN{{\Bbb N}}

\newcommand\proof{{\noindent\bf Proof.\ }}

\newtheorem{theorem}{Theorem}[section]
\newtheorem{lemma}[theorem]{Lemma}
\newtheorem{corollary}[theorem]{Corollary}
\newtheorem{proposition}[theorem]{Proposition}
\newtheorem{question}[theorem]{Question}
\newtheorem{re}[theorem]{Remark}

\newtheorem{definition}[theorem]{Definition}
\newtheorem{conjecture}[theorem]{Conjecture}

\newtheorem{notation}[theorem]{Notation}

\newenvironment{remark}{\begin{re}\em}{\end{re}}

\def\rank{{\rm rank\,}}
\def\d{{\rm d}}
\newcommand{\trans}[1]{\hspace{0,1cm} #1^{\hspace{-0,4cm}t}\hspace{0,24cm}}

\begin{document}
\title {Geometric stability of the cotangent bundle and the universal
cover of a
projective manifold} \author{Fr\'ed\'eric Campana and Thomas
Peternell\\ (with an appendix by Matei Toma)}
\date{November 8, 2007}

\maketitle

\vspace*{-0.5in}\section*{Introduction}
Let $X_n$ be a complex projective $n$-dimensional manifold and
$\tilde X$ its universal cover. The Shafarevich conjecture asserts
that $\tilde X$ is holomorphically convex, i.e. admits a
proper holomorphic map onto a Stein space. There are two extremal
cases, namely that this map is constant, i.e. $\tilde X$ is
compact. This means that $\pi_1(X)$ is finite and not much can be
said further. Or the map is a modification, i.e through the
general point of $\tilde X$ there is no positive-dimensional compact
subvariety, e.g. $\tilde X$ is Stein.
This happens in particular for $X$ an  Abelian variety or a quotient
of a bounded domain. It is conjectured (see [Ko93], and [CZ04] for
the K\"ahler case) that $X$ should then have a holomorphic submersion
onto a variety of general type with Abelian varieties as fibres,
after a suitable finite \'etale cover and birational modification.
This follows up to dimension $3$ from the solutions of the
conjectures of the Minimal Model Program.   We prove here a weaker
statement in every  dimension:

\begin{theorem} Let $X$ be a normal projective variety with at most rational
singularities.\\
(1) Suppose that the universal cover of $X$ is not covered by its
positive-dimensional
compact subvarieties. Then $X$ is of general type if
$\chi(\sO_X) \ne 0.$ \\
(2) In particular, if $X$ hat at most terminal singularities and $\tilde X$
is Stein (or does not
contain
a compact subvariety of positive dimension), then either
$K_X$ is ample, or we are in the following situation: $ K_X$ is
nef, $K_X^n = 0$, and $\chi(\sO_X) = 0$, with $n = \dim X$.
\end{theorem}

\noindent
This theorem is deduced (via the comparison theorem [Ca95], which
relates the geometric positivity of subsheaves in the cotangent
bundle to the
geometry of $\tilde X$) from the following more general:

\begin{theorem} Let $X$ be a projective manifold. Suppose that
$\Omega^p_X$ contains for some $p$ a subsheaf whose determinant is
big, i.e.
has maximal Kodaira dimension $n=dim(X)$. Then $K_X$ is big, i.e. $\kappa
(X) = n.$
\end{theorem}
\noindent
This in turn follows from a slight generalization of Miyaoka's
theorem that the cotangent bundle of a projective manifold is
``generically nef'' unless the manifold is uniruled, and from a
characterization
of pseudo-effective line bundles by moving curves
[BDPP04]. We show indeed that quotients of $\Omega_X^p$ have a
pseudo-effective determinant if $X$ is not uniruled and even more generally

\begin{theorem} Let $X$ be a projective manifold, $(\Omega^1_X)^{\otimes m}
\to \sS$ a torsion free
quotient. Then $\det \sS$ is pseudo-effective unless $X$ is uniruled.
\end{theorem} 

This
uniruledness criterion has also other applications, e.g. one can
prove that a
variety admitting a section in a tensor power of the tangent bundle
with a zero, must be uniruled. \\
Theorem 0.2 is actually a piece in a larger framework. To explain
this, we consider subsheaves $\sF \subset \Omega^p_X$
for some $p > 0.$ Then one can form $\kappa (\det \sF)$ and take the
supremum over all $\sF.$ This gives a
refined Kodaira dimension $\kappa^+(X)$, introduced in [Ca95]. Conjecturally
$$ \kappa^+(X) = \kappa (X) \eqno (*) $$
unless $X$ is uniruled. Theorem 0.2 is nothing but this conjecture
in case $\kappa^+(X) = \dim X.$ \\
We shall prove the conjecture (*) in several other cases.
It is actually a consequence of
the following more general:

\vskip .2cm \noindent {\bf Conjecture:}
{\it Suppose $X$ is a projective manifold, and suppose a decomposition
$$ NK_X = A + B$$
with some positive integer $N$, an effective divisor $A$ (one may assume $A$
spanned) and a pseudo-effective line bundle $B$.
Then $$\kappa (X) \geq \kappa (A).$$ }
\vskip .2cm \noindent

The special case $A=\sO_X$ implies that
$\kappa(X)\geq 0$ if $X$ is not uniruled, using the preceding result,
and the pseudo-effectiveness of $K_X$ when $X$ is not uniruled
([BDPP04]).

In another direction we obtain the special case in which $B$ is
numerically trivial:

\begin{theorem} Let $X$ be a projective
complex manifold, and $L\in {\rm Pic}(X)$ be numerically trivial.
Then:
\begin{enumerate}
\item $\kappa(X,K_X+ L)\leq \kappa(X)$.
\item If $\kappa(X)=0$,
and if $\kappa(X,K_X+ L)= \kappa(X)$, then $L$ is a torsion element
in the group ${\rm Pic}^0(X)$.
\end{enumerate}
\end{theorem}

\noindent In particular, if $mK_X$ is numerically equivalent to an
effective divisor, then $\kappa (X) \geq 0.$

\vskip .5cm \noindent
{\it Acknowledgement.} Our collaboration has been made possible by
the priority program {\it  ``Global methods in complex
geometry''} of the Deutsche Forschungsgemeinschaft, which we
gratefully acknowledge.

We thank P. Eyssidieux, T.Eckl, J.Stix for pointing out a gap in the first
version, and also C. Mourougane for his interesting observation
on our previous remark 3.6.

\vfill \eject

\tableofcontents

\section{Uniruledness Criteria}

\

Our main tool, of independent interest, is the generalisation \ref{unir}
below of Miyaoka's uniruledness criterion  \ref{miy}, which we first recall.

\

\begin{definition} Let $X$ be a complex projective $n$-dimensional manifold.
A vector
bundle $E$ over $X$ is
{\it generically nef}, if for all ample line bundles $H_1, \ldots,
H_{n-1}, $ for all $m_i $ sufficiently large
and for general curves $C$ cut out by $m_1 H_1, \ldots,
m_{n-1}H_{n-1}, $ the bundle $E \vert C$ is nef.
\end{definition}

\noindent Miyaoka's criterion [Mi87], with a short proof in [SB92], is the
following.

\begin{theorem}\label{miy} The cotangent bundle of a projective manifold
is generically nef if $X$ is not uniruled.
\end{theorem}

\

Before staing the first generalisation in \ref{unir} below, we need to
introduce the notion of movable class of curves, generalising complete
intersections curves.

We will denote by ${\overline {ME}}(X) $ the closed cone of (classes of)
movable curves, as defined in [BDPP04]. This is the smallest
closed cone containing all the classes of movable curves: a curve $C$ is
movable if it belongs to a covering family $(C_t)_{t \in T}$
of curves which is to say that $T$ is irreducible and projective, the
general $C_t$ is irreducible and the $C_t$ covers $X.$ \\

One of the main results of [BDPP04] is that ${\overline {ME}}(X)$ is the
closed convex cone generated by  classes $\alpha$ of the form $\alpha \pi_*(H_1 \cap \ldots \cap H_{n-1})$, with $\pi: X' \to X$ a modification
and $H_j$ very ample on $X',$ see (1.8) below. \\

Let $\alpha \in {\overline {ME}}(X).$ The slope of a torsion free sheaf
$\sE$ of rank $r$ with respect to $\alpha$ is defined by
$$ \mu_{\alpha}(\sE) = {{c_1(\sE) \cdot \alpha} \over {r}}.$$
A torsion free sheaf is $\alpha-$semi-stable, if for all proper non-zero
coherent subsheaves $\sF \subset \sE:$
$$  \mu_{\alpha}(\sF) \leq \mu_{\alpha}(\sE). $$

The general properties of $\alpha$-slopes are very much parallel to the
classical polarized case $\alpha = H_1 \cdot \ldots \cdot H_{n-1}$
with ample line bundles $H_i.$

\begin{proposition}\label{hn}
Let $X$ be a projective manifold and $\alpha \in {\overline {ME}}(X).$
Let $\sE$ be a non-zero coherent torsion free sheaf on $X$. Then:
\begin{enumerate}
\item  When $\sF$ ranges over all nonzero proper coherent subsheaves of
$\sE$, the slope $\mu_{\alpha}(\sF)$ is bounded from above.
Let $\mu_{\alpha}^{max}(\sE)$ be the maximum value.
\item There exists a unique largest subsheaf $\sE^{max}\subset \sE$ such
that 
$$\mu_{\alpha}(\sE^{max}) = \mu_{\alpha}^{max}(\sE).$$ The quotient
$\sE/\sE^{max}$ is torsion free.
\item Define inductively
$$\sE_0 =\lbrace 0\rbrace\subset \sE_1 = \sE^{max}\subset \ldots \subset
\sE_{s+1} = \sE$$ such that
$(\sE_{j+1}/\sE_j) = (\sE/\sE_j)^{max}$, for $j=0, \ldots ,s$.
This sequence is called the Harder-Narasimhan filtration of $\sE$ relative
to $\alpha.$ We write
$$\mu^{min}_{\alpha}(\sE):=\mu(\sE/\sE_s).$$
The quotients $\sE/\sE_j$ are the $\alpha$-semistable pieces of the
HN-filtration of $\sE$ relative to $\alpha$.
\item $\mu_{\alpha}(\sE_{j+1}/\sE_j) = \mu_{\alpha}^{max}(\sE/\sE_j) >
\mu_{\alpha}(\sE/\sE_{j+1})$, for $j \geq 0$.
\item $Hom(\sE_j,\sE/\sE_j) = 0 $ for all $j\geq 0$.

\item Let $\alpha \in  \overline{ME}(X)$ and $\sE$ and $\sF$ be
$\alpha-$semi-stable torsion free sheaves on $X$.
Then $\sE \hat \otimes \sF:= (\sE \otimes \sF)/{\rm tor}$ is again
$\alpha-$semi-stable.
\item $Hom(\wedge ^2\sE_j,\sE/\sE_j) = 0$ for all  $j\geq 0$.
\end{enumerate}
\end{proposition}

\proof The proof of the first four statements is essentially the same as in
the classical case of polarised varieties, see e.g. [MP97,p.42].
The last two properties follow from property $(4)$,
and the fact (see also [SB92]) that ${\rm Hom}(\sE,\sF) =0$ if
$\mu_{\alpha}^{min}(\sE) > \mu_{\alpha}^{max}(\sF)$. Property (6) is more
delicate, and proved in the appendix. For property (7), we proceed in the
usual way (see [SB92]), using (6).
\qed

The first generalization of Theorem 1.2 is

\begin{theorem}\label{unir} Let $X$ be a connected projective manifold, and
$\alpha \in {\overline {ME}}(X)$ of the form
$$\alpha = \pi_*(H_1 \cap \ldots \cap H_{n-1})$$
with $\pi: X' \to X$ a modification and $H_j$ very ample on $X'.$
If there exists a torsion free quotient sheaf
$$\Omega^1_X \to Q \to 0 $$ such that $c_1(Q) \cdot \alpha < 0$, then $X$ is
uniruled. \\
In other words, if $(C_t)$ is a covering family of curves which is the
birational image of hyperplane sections with $c_1(Q) \cdot C_t < 0,$ then
$X$ is uniruled. 
\end{theorem}

\begin{remark} {\rm
\begin{enumerate}
\item  The last assumption in Theorem 1.4  cannot be weakened to assuming
that, for generic $t\in T$, the bundle $\Omega^1_{X\vert C_t}$ is not nef
(i.e. $\Omega^1_X \vert C_t$ has a quotient $Q_t$ such that $\deg(Q_t)<0)$.
See [BDPP04], Theorem 7.7.
\item The last assumption is however satisfied if, for generic $t\in T$,
$\Omega^1_{X\vert C_t}$ is not nef, provided $C_t$ is an ample curve
obtained as 
intersection of $(n-1)$ generic members of a sufficiently high multiple of
some polarisation $H$ on $X$. This is a consequence of [MR82]. See [SB92].
\end{enumerate}}
\end{remark}

\begin{question} Let $X$ be a projective manifold and  $\pi:X'\to X$ be a
modification from another projective manifold $X'.$
Is $\pi^*(\Omega^1_{X})$ generically nef if $X$ is not uniruled?
\end{question}

The problem is to show that the last assumption of \ref{unir}  is satisfied,
if $C_t=\pi_*(C'_t)$, where $C'_t$ is a sufficiently ample curve on $X'$,
as in the preceding remark (2).

\vskip .2cm \noindent
\proof (of \ref{unir}) The proof follows the line of argumentation in
[SB92], using the notion of Harder-Narasimhan filtration for  $\alpha \in
{\overline {ME}}(X).$
Observe that we cannot use [MR82] in our context.
To be more precise, assume that $X$ is not uniruled. Then $K_X \cdot \alpha
\geq 0$ by [BDPP04], see Theorem 1.8 below.
Hence $\Omega^1_X $ is not $\alpha-$semi-stable and so is its dual $T_X.$ We
now define $\sF \subset T_X$ to be the largest piece $\sF_j$ of the
HN-filtration
of $T_X$ with $\mu(\sF_{j+1}/\sF_j) \geq 0,$ noticing that $c_1(\sF) \cdot
\alpha > 0.$ Let $\sG = T_X/\sF.$ Then Proposition 1.3(7) applies and
we conclude that $\sF$ is Lie closed. \\
As in [SB92] we now reduce to char $p$ and want to prove that $\sF$ is
$p-$closed. So
let $F: X \to X$  denote the absolute Frobenius,; we need to prove that
$${\rm Hom}(F^*(\sF),\sG) = 0.$$
Instead of restricting to curves as in [SB92] - which will not work in our
situation - we first observe that [SB98,Prop.1] remains true with exactly
the 
same proof in our situation. Then
we use the arguments of the first few lines of [La04, 2.5], substituting
[La04,2,4] by [SB98,Prop.1] to show that
$$\mu^{\rm max}_{\alpha}(F^*(\sF)) - \mu^{\min}_{\alpha}(F^*(\sF)) $$
is bounded independently of $p:$
$$ \mu^{\rm max}_{\alpha}(F^*(\sF)) - \mu^{\min}_{\alpha}(F^*(\sF)) \leq
({\rm rk}(\sF) - 1)H \cdot \alpha $$
for some fixed sufficiently ample line bundle $H$. This implies the
$p-$closedness of $\sF$, analogously to [SB92,9.1.3.5]
\vskip .2cm \noindent
Thus $\sF$ is Lie closed and $p-$closed, and therefore  $F$ is a 1-foliation
in the terminology of [SB92]. Hence one can form the quotient $\rho: X \to Y
= X/\sF$,
and obtain a $Y$-covering family of curves $C'_t = f_*(C_t)_{t\in T} $ plus
a line bundle $\sG'$ on $Y$ with $\det \sG
= \rho^*(G).$ Furthermore, we have
$$ C_t \cdot \rho^*(-K_Y) = p (C_t \cdot c_1(\sF)) + (C_t' \cdot
c_1(\sG)).$$
Then we proceed almost verbatim as in [SB92] to produce rational curves
through the general point of $X$ with a bound on their degree (with respect
to any polarisation on $X)$, which is independent of $p.$

\qed
We shall need the following generalization

\begin{theorem} \label{unir'}
 Let $X$ be a connected projective manifold, and  $\alpha \in {\overline
{ME}}(X)$ of the form
$$\alpha = \pi_*(H_1 \cap \ldots \cap H_{n-1})$$
with $\pi: X' \to X$ a modification and $H_j$ very ample on $X'.$
If there exists a torsion free quotient sheaf
$$(\Omega^1_X)^{\otimes m} \to Q \to 0 $$ such that $c_1(Q) \cdot \alpha <
0$, then $X$ is uniruled.
\end{theorem}

\proof As in the proof of Theorem 1.4, $(\Omega^1_X)^{\otimes m}$ is not
$\alpha-$semi-stable; let $\sS_m$ be the maximal destabilizing subsheaf.
>From our assumption
$$ \mu^{\rm max}_{\alpha}((\Omega^1_X)^{\otimes m}) = \mu_{\alpha}(\sS_m) >
\mu_{\alpha}((\Omega_X^1)^{\otimes m}) > 0.$$
Hence by Theorem 5.1 of the appendix,
$\Omega^1_X$ is not $\alpha-$semi-stable. Let $\sS_1 \subset \Omega^1_X$ be
the maximal  destabilizing subsheaf with torsion
free quotient $Q_1.$ By Corollary 5.4, we obtain
$$ \mu^{\rm max}_{\alpha}(\Omega_X^1) = \mu_{\alpha}(\sS_1) > 0.$$
Hence
$$ c_1(Q_1) \cdot \alpha < 0, $$
and $X$ is uniruled by Theorem 1.4.

\qed

\vskip .2cm  
\noindent Now we can strengthen the preceeding result, using
[BDPP04] (and answering a question asked in that paper). \\
First recall that a line
bundle $L$ on a projective manifold is
called {\it pseudo-effective} iff $c_1(L)$ is in the
closure of the cone generated by the (numerical equivalence classes
of the) effective divisors on $X.$

\noindent We will need the following result from [BDPP04] which will
also be crucial for Theorem 2.3.

\begin{theorem} Let $X_n$ be a projective manifold and $L$ a line
bundle on $X.$ Then $L$ is pseudo-effective if and only if
the following holds. Let $\pi: \hat X \to X$ be a birational map from
a projective manifold $X$. Let $H_1, \ldots, H_{n-1} $ be
very ample line bundles on $\hat X.$ Then
$$ L \cdot \pi_*(H_1 \cap  \ldots \cap H_{n-1}) \geq 0.$$
\end{theorem}

\noindent Together with theorem (1.7), this implies:

\begin{theorem} Let $X$ be a projective manifold and suppose that
$X$ is not uniruled. Let  $ Q$ be a torsion free quotient of
$(\Omega^1_X)^{\otimes m}$ for some $m>0$.
Then $\det Q$ is
pseudo-effective.

\end{theorem}

Here we made use of the following

\begin{notation} Let $\sF $ be a coherent sheaf
of rank $r$ on the
connected manifold $X$. We define its determinant -
a line bundle since $X$ is smooth -
to be
$$\det \sF = (\bigwedge^r \sF)^{**}.$$
\end{notation}

\proof
In order to show the pseudo-effectivity of $\det Q$, it suffices by
(1.4) to verify the following. Let $\pi:
\tilde X \to
X$ be birational from the projective manifold $\tilde X.$ Let $H_1,
\ldots, H_{n-1}$ be very ample on $\tilde X$ (and general in their
linear systems). Then
$$ \det Q \cdot \pi_*(H_1 \cap \ldots \cap H_{n-1}) \geq 0. \eqno (*)$$
So let us verify (*).
Eventually after replacing $H_i$ by large multiples and by setting
$$\tilde C = H_1 \cap \ldots \cap H_{n-1}, $$
Theorem 1.7  applies and $\pi^*\Omega^1_X \vert \tilde C$ is nef.
Hence $\pi^*(Q) \vert C$ is nef, therefore also $(\pi^*(Q)/{\rm torsion})
\vert \tilde C$ is nef, so that
$(\det \pi^*(Q))\vert C$ is nef. Denoting $r$ the (generic) rank of $Q,$ we
have $\bigwedge^r(\pi^*(Q)) = \pi^*(\bigwedge^r Q),$
hence we obtain a canonical map
$$ \bigwedge^r \pi^*(Q) \to \pi^*(\det Q).$$
This yields an inclusion $$\det \pi^*(Q) \subset \pi^*(\det Q).$$ Hence
$\pi^*(\det Q) \vert \tilde C$ is nef, too, and (*) is verified. This
finishes the proof.

\qed

\noindent Now a pseudo-effective line bundle is nef on moving curves; here
``moving'' means that the
deformations of the curve cover the variety.
Actually by [BDPP04] the closed cone
generated by by numerical equivalence classes of movable curves coincides
with the cone generated by classes of
``strongly movable'' curves. These are just the curves of the form $\pi(\hat
C)$, where  $\pi: \hat X \to X$ is
a modification, and $\hat C \subset \hat X $ is a generic intersection of
very
ample divisors $m_iH_i, 1 \leq i \leq n-1$ on $\hat X$. So we can state:

\begin{corollary} Let $X$ be a projective manifold and
suppose that $X$ is not uniruled. Let
$(C_t)_{t\in T}$ be an algebraic family of curves, parametrised by
the irreducible projective variety $T$. Assume this family is
covering (i.e.: the union of the $C_t$'s is $X$, and its generic member
is irreducible). \\
Let  $\sF$ be a torsion free quotient of
$(\Omega^1_X)^{\otimes m}$ for some $m>0$.
Then $c_1(\sF).C_t \geq
0$. \end{corollary}

\begin{corollary}  Let $X$ be a projective manifold and $L$ a
topologically trivial line bundle on $X$.
Let $m$ be a positive integer and
$$ v \in H^0(T_X^ {\otimes m} \otimes L) $$
a section with zeroes in codimension 1. Then $X$ is uniruled. \\
More generally, suppose that $\sF \subset T_X^{\otimes m} $
is a coherent subsheaf of rank $r$ such that $\det \sF$ is pseudo-effective
and that $\det \sF \to \bigwedge^r(T_X^{\otimes r})$
has zeroes in codimension 1. Then $X$ is uniruled.
\end{corollary} 

\begin{remark} {\rm
A classical result in group actions on a projective manifold $X$ says that
if
$X$ carries a holomorphic vector field with zeroes, then $X$ is uniruled.
If Question 1.6 had a positive answer, then we would be able to generalize
this result to arbitrary tensor
powers of the
tangent bundle, and we may also allow a twist with a topologically
trivial line bundle. In other words, we would be able to generalize (1.12)
by assuming there only the existence of some zero without saying
anything on the dimension of the zero locus. \\
\vskip .1cm \noindent
In fact, choose $p \in X$ such that $v(p) = 0.$ Let $\pi: \hat X \to X$
be the blow-up of $X$ at $p.$ Assume that $X$ is
not uniruled. Then, supposing that (1.6) has a positive answer,
$\pi^*(\Omega^1_X) $ is generically nef. Hence if
$\hat C$ is the curve cut out by sufficiently
general very ample divisors, then $\pi^*(\Omega^1_X) \vert \hat C $
is nef. Thus $\Omega^1_X \vert C$ is nef,
where $C = \pi(\hat C).$ Now $\hat C $ meets the exceptional divisor
of $\pi$ in a finite set, hence
$p \in C.$ In total, $(\Omega^1_X)^{\otimes m} \otimes L^* \vert C $
is nef, but its
dual has a section with zeroes.
This is impossible. So $X$ is uniruled.}
\end{remark}

\section{A characterization of varieties of general type}

\subsection{Refined Kodaira Dimension}

The following ``refined Kodaira dimension'' was introduced in
[Ca95]. It measures the geometric positivity of the cotangent bundle,
and not only  that of the canonical bundle. (Its definition is
justified in the next subsection).

\begin{definition} Let $X$ be a compact (or projective ) manifold.
Then $\kappa^+(X)$ is the maximal number
$\kappa (\det \sF) $, where $\sF \subset \Omega^p_X$ for $1 \leq p \leq
\dim X$ is a (saturated) coherent subsheaf.
\end{definition}

\noindent
Obviously we have $\kappa^+(X)\geq \kappa(X)$ for any $X$. \\
Assuming the standard conjectures of the Minimal Model Program, one
can easily describe
$\kappa^+(X)$ as follows (see [Ca95] for details,
where the following conjecture was
formulated):

\begin{conjecture}\label{conj} Let $X$ be a projective
manifold.
   If  $X$ is not uniruled (or if $\kappa(X)\geq 0$), then
$ \kappa^+(X) = \kappa (X).$
\end{conjecture}

\noindent
When $X$ is uniruled, one has
$$\kappa^+(X)=\kappa^+(R(X)),$$
where $R(X)$ is the so-called ``rational
quotient" of $X$; see [Ca95]. This rational quotient is not uniruled,
and so should be either one point or have
$\kappa^+(R(X))=\kappa(R(X))\geq 0$. Thus if $X$ is uniruled, one has
$\kappa(X) = - \infty$
but $\kappa^+(X) \geq 0$, {\it unless} $R(X)$ is one point, which
means that $X$ is rationally connected. In this latter case
$\kappa^+(X)=-\infty$. Conversely, if $\kappa^+(X) = - \infty,$ then $X$
should be 
rationally
connected. \\
Notice that $\chi(\sO _X)=1$ if
$\kappa^+(X)=-\infty$, because $h^0(X,\Omega^p_X)=0$ for $p>0$. In
[Ca95] it is shown that $X$ is simply connected if
$\kappa^+(X)=-\infty$ which of course is also true for $X$ rationally
connected. 
\vskip .2cm \noindent
The above conjecture is a
geometric version of the stability of the cotangent bundle of $X$
when $X$ is not uniruled. It is a version in which positivity of subsheaves
is measured by the Kodaira dimension of the determinant bundle, and
not  by the slope after restricting to ``strongly movable curves''.

\subsection{A Characterisation of Varieties of General Type}

As a consequence of the preceeding criteria for uniruledness, we first
solve the above conjecture in the extremal case when $\kappa^+(X)=n$
(we shall study in the next section below the intermediate cases):

\begin{theorem} Let $X_n$ be a projective manifold and suppose
$\kappa^+(X) = n$, i.e. some $\Omega_X^p$ contains a subsheaf
$\sF$ with $\kappa (\det \sF) = n.$ Then $\kappa (X) = n.$
\end{theorem}

\proof
First let us see that $X$ is not uniruled. In fact, otherwise take
a covering family of rational curves and select a
general member $C$ so that $T_X \vert C $ is nef. Hence the dual of
$\Omega^p_X \vert C$ is nef and therefore $\sF \vert C$
cannot have ample determinant. So $X$ cannot be uniruled.\\
Of course, we may assume that $\sF$ saturated, hence $Q = \Omega^p_X / \sF$
is
torsion free. By taking determinants we get
$$ mK_X = \det \sF + \det Q$$
for some positive integer $m.$
We learn from (1.6) above that $\det Q$ is pseudo-effective. Thus $K_X$ is
big, as a sum of a big and a pseudo-effective divisor.\\
\qed

\subsection{The intermediate case}

\noindent In this section we want to study the above conjecture \ref{conj}
in
the intermediate case $n>\kappa(X_n)\geq 0$. \\
We shall reduce
Conjecture \ref{conj} to (special cases of) a seemingly considerably
simpler:

\begin{conjecture} \label{k=a+b} Let $X$ be a projective
manifold. Let $N K_X = A + B$ with some positive integer $N>0$, $A$
effective 
and $B$ pseudo-effective. Then $\kappa (X) \geq \kappa
(A).$
\end{conjecture}

\begin{remark} (1) By suitably blowing up, it is easily seen that
Conjecture 2.4 is equivalent to the analogous conjecture
with $A$ {\it always assumed to be spanned.}
\vskip .2cm \noindent
(2) If $\nu(L)$ denotes the numerical
dimension of an arbitrary
pseudo-effective line bundle as introduced by Boucksom [Bo02], then
the generalised abundance conjecture states
$$ \kappa (K_X) = \nu(K_X).$$
If this generalised abundance conjecture holds, then Conjecture
\ref{k=a+b} holds when $\kappa(X)=0$,
a case sufficient to imply
conjecture \ref{conj} (see below). In fact, if $\kappa (K_X) = 0$
and $N K_X = A + B$ with $A$ spanned and $B$ pseudo-effective, then
$\nu(A+B) = 0$, hence $\nu (A) = 0$ and therefore $A = 0, A$ being spanned.
\end{remark}

\noindent We start with an immediate observation:

\begin{proposition} Conjecture \ref{k=a+b} implies Conjecture
\ref{conj}, when $X$ is not uniruled (and so when $\kappa(X)\geq 0$ ).
\end{proposition}

\proof Let
$\sF$ be a saturated subsheaf of $\Omega ^p_X$ such that $\kappa(X,
\det(\sF))=\kappa^+(X)\geq 0$, then  $Q = \Omega ^p_X / \sF$ is
torsion free. By taking determinants we get
$$ mK_X = \det \sF + \det Q$$
for some positive integer $m.$
We know that $\det Q$ is pseudo-effective, because $X$ is not
uniruled. By Conjecture \ref{k=a+b}, we get the claim, since
$A := \det(\sF)$ is $\bQ$-effective.
\qed

\noindent We now show that Conjecture \ref{k=a+b} (in case
$\kappa(X)\geq 0$), and so
\ref{conj}, is implied by the special case $\kappa(X)=0$ of Conjecture
\ref{k=a+b}. More precisely:

\begin{proposition} \label{k^+}Let
$X_n$ be a projective manifold with $\kappa (X)
\geq 0$. Let $d = n-\kappa(X) \geq 0.$ If  Conjecture \ref{k=a+b}
holds for all
manifolds  $G$ of dimension $d$ and with
$\kappa(G)=0$, then Conjecture \ref{k=a+b} (and so also Conjecture
\ref{conj}) holds for
$X$.
\end{proposition}

\proof
By blowing up we may assume that the Iitaka fibration
$g: X \to W$ is holomorphic. Let $G$ be a general fiber of $g$. Thus
$\kappa(G)=0$.
Let $A$ be effective and $B$ pseudo-effective on $X$ such that
$$ NK_X = A + B$$
for some positive integer $N.$
Then $A_G$ is effective,$B_G$ is pseudo-effective and
$$ NK_G = A_G + B_G.$$
Thus by Conjecture \ref{k=a+b} applied to $G$, we conclude that
$\kappa(G,A_{\vert G}\leq 0$.
By the easy additivity theorem for
Kodaira dimension, we obtain that
$$\kappa(X,A)\leq 
\dim(W)+\kappa(G,A_{\vert G}) \leq \dim(W)=\kappa(X).$$
\qed

\noindent The preceeding observation shows that the only two crucial cases
of 
Conjecture \ref{k=a+b} are when $\kappa(X)$ is either $0$, or
$-\infty$.

\noindent
We now give some cases in which Conjecture \ref{k=a+b} can be
solved, so that \ref{k^+} can be applied. \\
We first recall
a notion from Mori theory. Let $X$ be a
projective manifold. A variety $X'$ with at most terminal
singularities is said to be a {\it good minimal model} for $X$, if
$X'$ is birational to $X$ and some $mK_{X'} $ is (locally free and)
spanned. Good minimal models are predicted to exist for every $X$
with $\kappa (X) \geq 0$ but this known only in dimension
up to $3.$

\begin{proposition} Let $G$ be a projective manifold with $\kappa (G)
= 0.$ Suppose $G$ has a good minimal model and that
$$ N K_G = A + B$$
with $A$ effective and $B$ pseudo-effective.
Then $\kappa(A) = 0.$
\end{proposition}

\proof Let $G'$ be a good minimal model for $G.$ Then $K_{G'} \equiv
0$ and actually $K_{G'}$ is torsion.
Choose a smooth model $\hat G$ with holomorphic maps $\pi: \hat G \to G$ and
$\lambda: \hat G \to G'.$ There is an effective divisor $E$ supported on the
exceptional locus of
$\pi$ such that $K_{\hat G} = \pi^*(K_G) + E.$
Then we can write 
$$ NK_{\hat G} = \hat A + \hat B$$
with $\hat A = \pi^*(A) + NE$ effective and $\hat B = \pi^*(B) $
pseudo-effective. 
Now consider $A' = \lambda_*(\hat A)$ and $B' = \lambda_*(\hat B).$ Then
$A'$ is effective, $B'$ is
pseudo-effective and
$$NK_{G'} = A'+B'.$$
It follows $A' = B' = 0$ so that $\kappa(A) = 0.$
\qed

\noindent Since good minimal models exist in dimension up to 3, Prop.
\ref{k^+}
gives in particular:

\begin{theorem} Let $X_n$ be a projective manifold, $\kappa (X)
\geq 0.$ Suppose $\kappa(X) \geq n-3.$
Then $\kappa^+(X) = \kappa (X).$ \end{theorem}

\vskip .2cm \noindent
For some other result towards
(2.4) we state

\begin{proposition} Let $X_n$ be a projective manifold, $N K_X = A +
B$ with $A$ spanned and $B$ pseudo-effective.
Let $f: X \to Y$ be the fibration determined by $\vert A \vert.$ Let
$F$ be the general fiber of $f.$ If $B \vert F$ is
big, then $K_X$ is big, i.e. $\kappa (X) = n.$
\end{proposition}

\proof This is proved in [CCP05].
\qed

\begin{corollary} Let $X_n$ be a projective manifold, $N K_X = A + B$
with $A$ spanned and $B$ pseudo-effective. If $\kappa (A) n-1,$ then $\kappa (X) \geq n-1.$
\end{corollary}

\proof Let $f: X \to Y$ be the fibration associated with $A$ and let
$F$ denote the general fiber. Since $\dim F = 1,$ either
$B_F$ is ample or $B_F \equiv 0.$ \\
In the first case we simply apply (2.10). In the second we notice $N
K_F = B_F \equiv 0$ so that $F$ is elliptic and $B_F = 0.$ Then
we can write
$$ mB = f^*(L) + \sum d_i D_i $$
with $L$ a line bundle on $Y,$ with $d_i$ integers, not necessarily
positive, and with $D_i$ irreducible divisors with
$\dim f(D_i) \leq n-2,$ but not pull-backs of divisors on $Y.$
Intersecting with movable curves in $Y$, it is easy to see that $L$
is pseudo-effective. Writing $A = f^*(A'),$ it follows that
$A' + L$ is big. Hence
$$ \kappa(NmK_X + \sum (-d_i'D_i)) = n-1, $$
where $d_i'$ are just the negative $d_i.$ Then however
$$ \kappa (X) = n-1, $$ too.
\qed

\section{Numerical properties of the Kodaira dimension}

We solve here Conjecture \ref{k=a+b} in the special case where $B$
is numerically trivial.

\begin{theorem} Let $X$ be a projective
complex manifold, and $L\in {\rm Pic}^0(X)$ be numerically trivial.
Then:
\begin{enumerate}
\item $\kappa(X,mK_X+ L)\leq \kappa(X)$.
\item If $\kappa(X)=0$,
and if $\kappa(X,mK_X+ L)= \kappa(X)$, then $L$ is a torsion element
in the group ${\rm Pic}^0(X)$.
\end{enumerate}
\end{theorem}

\begin{remark} The
conclusion of (2) above does no longer
hold when $\kappa(X)\geq 1$, as shown by curves (or even arbitrary
manifolds) of general
type.
\noindent
Another point not shown by our arguments is the behaviour of
the modified plurigenera
$p^+_m(X):={\rm sup} \lbrace h^0(X,mK_X+L),
L\equiv 0 \rbrace$, as $m$ is large and divisible. One may expect
that then $p^+_m(X)=p_m(X)$, and that the maximum is attained at a
torsion point, for every $m>0$ (this is true for $m=1$, by the
arguments below).
\end{remark}

\proof We proceed in two
steps:
\vskip .2cm \noindent
A. We prove the result in the special case where
$\kappa(X)\leq 0$. This will be done below in the two next
propositions.
\vskip .2cm \noindent
B. We now reduce the general case where $\kappa(X) >
0$ to the special case $\kappa(X) = 0$, as in \ref{k^+} above.

\vskip .2cm \noindent
Observe first that the statements involved are preserved by
birational transformations of $X$. We can thus assume that both $f,
g$ are holomorphic, where $g:X\to W$ is the Iitaka-Moishezon
fibration of $X$ defined by some $\vert mK_X\vert $, and  $f:X\to Y$
is the Iitaka fibration defined by some $\vert m(K_X+L) \vert $. If
$G$ is a general fibre of $G$, then it is sufficient to show that
$f(G)$ is a single point of $Y$. But then $f_{\vert G}$ is nothing,
but the Iitaka fibration on $G$ defined by $(K_X+L)_{\vert G}$.
Because $\kappa(G)=0$, we obtain the conclusion from the first step.

{\begin{remark} We see moreover that, in order to have equality
$\kappa(K_X+L)=\kappa(X)$, it is necessary that $L_{\vert G}$ be
torsion. This is Step B for the claim (2).
\end{remark}

\noindent To conclude the proof of the preceding
theorem, we still need to solve the case $\kappa(X)\leq 0$. This is
the content of the next two propositions (in which additive and
multiplicative notations for line bundles are mixed ).

\noindent We first deal
with the case $m=1$.

\begin{proposition} Let $X$ be a projective manifold, $L \in {\rm
Pic}^0(X).$ Suppose that $h^0(K_X \otimes L) \geq r $ for $r = 1$ or
$r = 2.$ Then
\begin{enumerate}
\item There exists a finite \'etale abelian cover
$\tilde X \to X$ such that $h^0(2K_{\tilde X}) \geq r.$ (In
particular, $\kappa(X)\geq r-1\geq 0$)
\item If $\kappa(X)=0$, and if
$h^0(K_X \otimes L) = 1$,  $L$ is a torsion element in the group
${\rm Pic}^0(X)$.
\end{enumerate}
\end{proposition}

\proof Fix $r = 1$ or $r = 2$. Our aim is to prove that there is a
finite \'etale cover $\tilde X \to X $ such that
$$ h^0(K_{\tilde X} \otimes L^*) \geq r, $$
Then the canonical morphism
$$ H^0(K_{\tilde X} \otimes L) \otimes H^0(K_{\tilde X} \otimes L^*)
\to H^0(2K_{\tilde X}) $$
will give the first claim $h^0(2K_{\tilde X}) \geq  r.$ \\
Let
$$ S = \{ G \in {\rm Pic}^0(X) \  \vert \ h^0(K_X \otimes G) \geq r \}. $$
If $S = {\rm Pic}^0(X),$ then we are clearly done. Otherwise Simpson
[Si93] gives the structure of $S:$
$$ S = \bigcup \{A_i +  T_i \}$$
with $A_i$  {\it torsion} elements and $T_i$ subtori of ${\rm Pic}^0(X).$ 
Let
$$ S^* = \{ G^* \ \vert \ G \in S \}.$$
Then $S^* = \bigcup \{-A_i+ T_i\}.$ Hence we can write $L^* = -A _i+
\tilde L $ with some $ \tilde L \in  T_i.$
Consequently $\tilde L = \hat L - A _j$ with $\hat L \in S $,  so that in 
total
$$ L^* = -A _i-A_j+ \hat L.$$
Now choose a finite \'etale cover $f: \tilde X \to X$ such that
$f^*(A_i+A_j) = \sO. $ Then we conclude
$$ h^0(K_{\tilde X} \otimes f^*(L^*)) \geq r $$
and we are done for assertion 1. \\
(Notice that if $L$ is unitary flat, then by Hodge theory it is
obvious that $h^0(K_X \otimes L) = h^0(K_X \otimes L^*)$, without
using  [Si93]).
\vskip .2cm \noindent
Let us now prove statement (2). Replacing $X$ by
$\tilde X$ as above, we can write: $L^*=\tilde L$, with
$\tilde
L\in T_i\subset S$. If $T_i$ is not trivial, we get a one-parameter
subgroup $L_t, t \in \bC$ contained in $T_i \subset S$. The canonical
morphisms
$$ H^0(K_{\tilde X} \otimes L_t) \otimes H^0(K_{\tilde X} \otimes
L_t^*) \to H^0(2K_{\tilde X}) $$
show that $h^0(2K_{\tilde X}) \geq  2$, contradicting our
assumption that $\kappa(X)=0$.\\

\qed

We shall now reduce the general case of $m \geq 2$ to the special
case $m=1$, by means of cyclic covers.\\

\begin{theorem} Let $X$ be a projective manifold and $L$ a line
bundle with $c_1(L) = 0$ in $H^2(X,\bZ).$
\begin{enumerate}
\item Suppose that there is a positive integer $m$ such that
$h^0(mK_X \otimes L) \geq 2.$ Then
$\kappa (X) \geq 1.$
\item Suppose that there is a positive integer $m$ such that
$h^0(mK_X \otimes L) \ne 0. $ Then
$\kappa (X) \geq 0.$
\item Suppose that $\kappa(X)=0$, and that
$h^0(mK_X \otimes L) \ne 0. $ Then $L$ is torsion in $\rm{Pic}^0(X).$
\end{enumerate}
\end{theorem}

\proof We only prove (1), the proof of (2) being identical,
simply omitting the divisor $D$ in the arguments below. Since our
claim is invariant by finite \'etale covers,
we can pass to such covers as we like. In particular, we may assume
that $L \in {\rm Pic}^0(X).$ If $m = 1,$ then our claim is Proposition
4.1, hence we shall assume $m \geq 2.$
Furthermore we may assume that
$L = mL', $ so that
$$ h^0(m(K_X \otimes L')) \geq 2. $$
Let $\sum b_i B_i $ be the fixed part of $\vert m(K_X \otimes L')
\vert,$ so that we can write
$$ m(K_X \otimes L') = \sum b_i B_i + D $$
with $D$ reduced and movable.
By possibly blowing up we may assume that the support of $\sum b_i
B_i + D$ has normal crossings. Now take the $m-$th root, normalize and
desingularize to obtain $f: Y \to X.$ We have to compute $f_*(K_Y),$
following [Es82,Vi83].
In fact, in

troduce the line bundles
$$H_j = j(K_X \otimes L') - \sum [jb_i m^{-1}] B_i. $$
Here [x] denotes the integral part of $x.$
Then
$$ f_*(K_Y) = K_X \otimes \bigoplus_{j = 0}^{m-1} H_j.$$
Hence the direct summand of $f_*(K_Y) \otimes L'$ corresponding to $j = m-1$
is just
$$ D + \sum_i (b_i - [b_i(m-1) m^{-1}]) B_i.$$
Since $D$ moves, we obtain
$$ h^0(f_*(K_Y) \otimes L') \geq 2,$$
hence $$ h^0(K_Y \otimes f^*(L')) \geq 2$$
so that $\kappa (Y) \geq 1.$
We still need to prove $\kappa (X) \geq 1.$
In to order to do that we must look more carefully at $f: Y \to X.$
This map decomposes as follows. First we take the cyclic covering
$h_0: Y_0 \to X$ determined by $m(K_X \otimes L') = \sO_X(D).$
Then we have the normalisation $h_1: Y_1 \to Y_0,$ and finally $h_2:
Y \to Y_1$ is a desingularisation.
Then $Y_0$ is Gorenstein and
$$ K_{Y_0} = h_0^*(mK_X \otimes (m-1)L');$$
furthermore $$ K_{Y_1} \subset h^*(K_{Y_0}) $$
via the trace map $(h_1)_*(K_{Y_1}) \to \sI \otimes K_{Y_0} $ (with
$\sI$ the conductor ideal) and finally
$$ (h_2)_*(K_Y) = K_{Y_1}  \eqno (+)$$
since $Y_1$ has rational singularities [Es82,Vi83].
In total $$ K_Y \subset f^*(mK_X \otimes (m-1)L')) + \sum a_i E_i$$
where $E_i$ are the exceptional components for $h_2$ and $a_i$ are
integers and we have
$$ (h_2)_*(\sO_Y(\sum a_i E_i)) = \sO_{Y_1} $$
by (+).
Thus
$$ K_Y \otimes f^*((1-m)L') \subset  f^*(mK_X) + \sum a_i E_i$$
We claim that there is a finite \'etale cover $g: \tilde Y \to Y$ such that
$$ h^0(K_{\tilde Y} \otimes g^*f^*((1-m)L')) = h^0(g^*(K_Y) \otimes
g^*f^*((1-m)L')) \geq 2. \eqno (*)$$
This will prove $\kappa (g^*(f^*(K_X) + \sum a_i E_i)) \geq 1,$ hence
$\kappa (X) \geq 1.$ \\
Now (*) is an easy application of Simpson's theorem, this time on
$Y.$ In fact, we introduce $S_Y$, the analogue of $S$ on $X.$
Then $S_Y = \bigcup_j \{B_j + T'_j\}$ with $B_j$ torsion and $T'_j$
tori. We already know that $f^*(L) \in S_Y$ so that
$$ f^*(L) = B + M $$
with $M \in T'_j$, and $B=B_j$  for some $j.$ Then $M' = (1-m)f^*(L)
\in T_j'$ and we find a finite \'etale cover $g: \tilde Y \to Y$
such that $g^*((1-m)f^*(L)) = g^*(M'),$ i.e. we choose $g$ such that
$g^*(B) = \sO.$ Now
$$ h^0(K_Y+B+M') \geq 2$$
since $M' \in T'_j.$ Hence $h^0(K_{\tilde Y}+g^*(M')) \geq 2$ which gives 
(*).
\vskip .2cm \noindent
As already said, the proof of (2) is essentially the same, omitting the
movable part.
\vskip .2cm \noindent
The proof of assertion (3) is the same as that of
assertion (2) of the preceeding Proposition 3.4.

\qed

\begin{remark} The preceeding result makes plausible the
expectation that the generalised Green-Lazarsfeld sets
$$S_{m,p,r}=\{L\in {\rm Pic}^0(X) \ \vert \  h^p(mK_X\otimes L)\geq r\}$$
might have the same structure as in [Si93] (finite union of translates of
subtori by torsion elements).

In fact, up to the word `` torsion" above, this is a consequence of the 
Abundance Conjecture, 
as C. Mourougane observed. Indeed he showed in [Mo99], thm. 5.3, that the 
Green-Lazarsfeld cohomological loci have this structure for ``good" 
divisors.
\end{remark}

\begin{corollary} Let $X$ be a projective manifold, $A$ effective and
$B$ pseudo-effective divisors on $X$. Assume that $mK_X = A + B$ for
some positive integer $m.$ Suppose also that $\nu (B) = 0$, in the
sense of [Bo02].  Then $\kappa (X) \geq 0.$
\end{corollary}

\proof By [Bo02], we can write $B \equiv \sum b_i B_i$ with positive
rational numbers $b_i.$ Now apply (3.1).

\qed

\section{The Universal Cover}

Another invariant of $X$ is defined via the universal cover
$\tilde X$ of a compact K\"ahler or projective manifold $X.$
By identifying points in $\tilde X$ which can be joined
by a compact connected analytic set, one obtains an almost holomorphic
meromorphic
map $\tilde X \rightharpoonup \Gamma (\tilde X) $. Here ``almost
holomorphic'' is to say that the degeneracy locus does not project
onto the image. If $\tilde X$ is holomorphically convex (which is
expected to be always true by the so-called
Shafarevitch conjecture), then this map is holomorphic and is just the usual
Remmert holomorphic reduction. In any case it induces the so-called
Shafarevich
map $\gamma_X: X \rightharpoonup \Gamma (X) = \Gamma (\tilde X) /
\pi_1(\tilde X).$

\begin{definition} $\gamma d(X) = \dim \Gamma (X) $ is the
$\Gamma-$dimension of $X$.
\end{definition}

\noindent Notice that $\gamma d(X) = 0 $ iff $\pi_1(X)$ is finite  and that
$\gamma d(X) = \dim X$ iff through the general point of $\tilde X$
there is no positive dimensional compact subvariety, i.e. $\tilde X$
geometrically seems as a modification of a  Stein space. \\
The following result [Ca95,(4.1)] gives a relation between
$\kappa^+(X)$ and $\gamma d(X).$

\begin{theorem} Let $X$ be a compact K\"ahler manifold. If
$\chi(X,\sO_X) \ne 0,$ then either
\begin{enumerate}
\item $ \kappa^+(X) \geq \gamma
d(X)$, or
\item $\kappa^+(X) = - \infty $, and so $X$ is simply
connected.
\end{enumerate}
\end{theorem}

\noindent By (2.9) we then obtain

\begin{corollary} Let $X_n$ be a projective manifold. Suppose that  $\kappa 
(X)
\geq n-3 $ and $\chi(\sO_X) \ne 0.$
Then $\kappa (X) = \kappa^+(X) \geq \gamma d(X).$
\end{corollary}

\noindent In particular, if $n = 4$, $\kappa (X) \geq 1$,
$\pi_1(X)$ is infinite
and $\chi(\sO_X) \ne 0,$ then $\kappa (X) \geq \gamma d(X)\geq 1.$
In other words, if $X$ is a projective 4-fold with $\kappa (X) = 0$
and $\pi_1(X) $ is not finite, then either $\chi(\sO_X) = 0$;
so there is either a holomorphic 1-form, or a holomorphic 
3-form, or: $\kappa^+(X)\in\{1,2,3\}$.

\noindent Hence as in [Ca95, 5.9], we conclude:

\begin{corollary} 
Let $X$ be a projective manifold of dimension $4$ such that 
$\kappa(X)=0$, and $\chi(\sO_X) \neq 0$. Then either $\pi_1(X)$ is finite 
and has at most $8$ elements, or $\kappa^+(X)\in\{1,2,3\}$.
\end{corollary}

\noindent This result should hold in arbitrary dimension $n$, with $8$ 
replaced by $2^{n-1}$, as 
a consequence of the standard conjecture that $\pi_1(X)$ should be 
almost abelian if $\kappa(X)=0$.

\noindent From Theorem 2.3 we deduce

\begin{corollary} Let $X$ be a normal projective  variety with at most 
rational singularities and suppose that
its universal cover is not covered by its positive-dimensional
compact subvarieties. Then  $X$ is of general type if $\chi(\sO_X) \ne 0.$
\end{corollary}

\proof If $X$ is smooth, then
by our assumption and (4.2), we have $\kappa^+(X) = \dim
X$ or $\chi(\sO_X) = 0.$ Now  theorem (2.3) gives the claim. \\
So it remains to reduce the general case to the smooth. Note that $\tilde X$ 
is irreducible since $X$ is normal. 
Consider a projective desingularisation $\pi: Y \to X$ and
let $\tilde \pi: \tilde Y \to \tilde X$ be the induced maps on the level of 
universal covers. Then $\tilde \pi$ is onto with discrete
fibers over the smooth locus of $\tilde X.$ Hence $\tilde Y$ is not covered 
by positive-dimensional compact subvarieties, too, because their
$\tilde \pi-$images would again be compact. By the solution of the smooth 
case, we either have $\chi(X,\sO_Y) = 0$ - hence $\chi(\sO_X) = 0$ by
the rationality of the singularities of $X$ - or $Y$, hence $X,$ is of 
general type.  
\qed

\begin{corollary} Let $X_n$ be a projective manifold or a normal projective 
variety with at most terminal singularities whose universal
cover is Stein (or has no positive-dimensional subvariety).
Then either $K_X$ is ample or
$\chi(\sO_X) = 0, $ $K_X$ is nef and $K_X^n = 0.$
\end{corollary}

\proof This is immediate from (4.5) by observing that $X$ does not
have any rational curve, so that $K_X$ must be nef by Mori
theory. Moreover if $K_X$ is big, then $K_X$ is ample by Kawamata [Ka92].

\qed

\noindent We are lead to ask for the structure of projective manifolds $X_n$
whose universal cover is Stein and with $K_X^n = 0$.

\begin{conjecture}\label{covst} Let $X_n$ be a projective manifold
whose universal
cover $\tilde X$ is Stein. Assume $K_X^n = 0.$
Then up to finite \'etale cover of $X,$ the manifold $X$ has a torus
submersion over a projective manifold $Y$ with $K_Y$ ample
and universal cover again Stein.
\end{conjecture}

\noindent If the universal cover is only assumed not to admit a
positive-dimensional subvariety through the
general point, then one expects a birational version of \ref{covst},
which is actually proved in
[Ko93,5.8]. Here is the ``Stein version'' of this result which does
not follow immediately from
Koll\'ar's result since we make a biholomorphic statement.  The main
point is to explain that we must
have a holomorphic Iitaka fibration which is ``almost smooth'' and
then apply Koll\'ar's techniques to make it smooth.

\begin{proposition} Conjecture \ref{covst} holds if $\kappa (X) \geq n-3$.
\end{proposition}

\proof
(1) Since the case $\kappa (X) = n-1$ is the simplest, we do
it first. Here the numerical dimension $\nu(X) = \kappa (X),$
so that $K_X$ is good, i.e. some multiple is spanned [Ka85b].
Therefore we have a holomorphic Iitaka fibration $f: X \to Y.$ The
general fiber is an elliptic curve. Since $X$ does not contain
rational curves, it follows easily that all fibers are elliptic,
sometimes multiple. Now [Ko93,sect.6] yields a
finite \'etale cover
such that the induced map is smooth; see below for some details. \\
(2) In the other case we consider the normalized graph $p: \sC \to X$
of the family determined by the general fibers of the
meromorphic Iitaka fibration. Let $q: \sC \to T$ denote the parameter space.
All irreducible fibers of $q$ have dimension 2 (resp. 3) and every
such fiber is an
\'etale quotient of a torus by Lemma 4.9 below.
Now we have a formula (via the trace map)
$$ K_{\sC} = p^*(K_X) + E $$
with an effective (Weil) divisor $E.$
Restricting to a general (normal, hence smooth by (4.8) below) fiber
$F$ of $q,$ we get
$$ 0 \equiv p^*(K_X) \vert F + E \vert F.$$
Hence $p^*(K_X) \vert F  \equiv  0  = E \vert F.$
Now consider the reduction $F_0$ of a component of a singular fiber
(or rather its normalization)
and use the conservation law (and the nefness of $p^*(K_X)$) to
deduce  $p^*(K_X) \vert F_0 \equiv 0$.
Thus $p^*(K_X)$ is
``$q-$numerically trivial''.
This proves immediately $\nu(X) = n-2$ (resp. $\nu(X) = n-3$) and
again $mK_X$ is
spanned for a suitable $m$.\\
Now let again $F_0$ be the reduction of a component of a singular
fiber $F$, this time
of the Iitaka fibration $f: X \to Y.$ \\
We claim that actually $F aF_0$ and that $f$ is equidimensional. \\
If $\dim T = 2$, this is easy and well-known
of course (take a general curve through $f(F_0)$ and observe that
singular non-multiple fibers produce rational curves). \\
So suppose $\dim T = 1.$ Take $\mu$ maximal such that
$\mu F_0 \subset F.$ Then $N_{F_0}^{*\mu}$ has a section which has a
zero, since $F$ is reducible. Hence $K_{F_0}
= - D$ with $D$ a $\bQ$-effective divisor by the adjunction formula.
Now normalize and then desingularize. The result
$\hat F_0$ has $\kappa (\hat F_0) ) = - \infty$ (use formula (*)
below), so that $F_0$ is uniruled. Since this is forbidden
by the universal cover, we obtain $F = aF_0.$ \\
Then $K_{F_0} \equiv 0,$ so that its normalization $\tilde F_0$ has
$$ K_{\tilde F_0} = \nu^*(K_{F_0}) - \tilde N \eqno (*)$$
with $\tilde N$ the preimage of the non-normal locus. Since
$K_{\tilde F_0} \equiv 0$ by (2.9), we conclude that
$F_0$ must have been normal, hence smooth. \\
Now we apply [Ko93,5.8] to obtain a finite \'etale cover $X'$ of $X$ which
is birational to a torus submersion. But since $X'$ does not contain 
rational
curves, we obtain a holomorphic birational map from a torus submersion to 
$X'.$
Since multiple fibers cannot be resolved by birational 
transformations on the base,
we conclude that $X'$ is a torus submersion itself.
\qed

It remains to prove the following lemma of independent interest.

\begin{lemma} Let $X$ be an irreducible reduced variety of dimension
at most 3. Assume that the universal cover of $X$ is Stein
(or does not contain compact subvarieties).
Let $\tilde X \to X$ be the normalization and $\pi: \hat X \to \tilde
X$ be a desingularization. Suppose $\kappa (\hat X) = 0. $
Then $\tilde X$ is an \'etale quotient of a torus.
\end{lemma}

\proof We only treat the case $\dim X = 3,$ the surface case being
easier and left to the reader.
By [NS95], $\hat X$ admits a finite
\'etale cover $h:X' \to \hat X$ which is birational to a product of a
simply connected manifold and an abelian variety. By our assumption
on the universal cover, the simply connected part does not appear.
It follows that the Albanese map $ \alpha: X' \to A$ is birational.
Now all irreducible components of all non-trivial fibers $\alpha$ are 
filled up by rational
curves ($\alpha$ factors via Mori contractions). Since $\tilde X$ 
does not contain rational curves, the map $X' \to \hat X \to \tilde X$
therefore factors over $\alpha,$ i.e. we obtain a finite map $g:A \to 
\tilde X.$ \\
This map is \'etale in codimension 1. In fact otherwise by the 
ramification formula $K_A = g^*(K_{\tilde X}) + R$ (as Weil divisors).
Thus $-K_{\tilde X}$ is non-zero effective and therefore $\kappa 
(\hat X) = - \infty,$ contradiction. \\
We want to see that $\tilde X$
is actually smooth and an \'etale quotient of $A$.
First notice that $\tilde X$ is $\bQ-$Gorenstein (if $g$ has degree 
$d$, then $dK_{\tilde X} = \sO $ on the regular part of $\tilde X,$
hence everywhere).
Now we can
compare the formulas
$$ K_{X'} = h^* \pi^*(K_{\tilde X}) + \sum a_i E_i $$
and
$$ K_{X'} = \sum b_j F_j $$
where $E_i$ are the preimages of the $\pi-$exceptional components and
$F_j$ are the $\alpha-$exceptional components; notice $b_j > 0.$
Then both sets of exceptional divisors are equal, and thus all $a_i >
0.$ Therefore $\tilde X$ has only terminal singularities. We also
notice that
$\pi_1(\tilde X)$ is almost abelian, i.e. abelian up to finite index.
Therefore $\pi_1(\tilde X)$ is abelian after finite \'etale cover.
Then [Ka85] applies and $\tilde X $ is an \'etale quotient of an
abelian threefold. Here of course we use again that the universal
cover of $\tilde X$ is Stein.
\qed

\section{Stability and tensor products} 

Recall that $\overline{ME}(X)$ denotes the movable cone of the 
$n-$dimensional projective manifold $X$. We say that $\alpha \in  
\overline{ME}(X)$
is {\it geometric}, if
there exists a modification $\pi: \tilde X \to X$ from the projective 
manifold $\tilde X$ and  ample line bundles $H_i$ such that
$$ \alpha = \lambda \pi_*(H_1 \cap \ldots \cap H_{n-1})$$ 
with a positive multiple $\lambda.$ 
By definition,  $\overline{ME}(X)$ is the closed cone generated by the 
geometric classes. \\ 
If $\sE$ and $\sF$ are torsion free sheaves,
then we put $\sE \hat \otimes \sF = (\sE \otimes \sF)/{\rm tor}.$ 
The first main result is well-known in case of an ample polarization. 

\begin{theorem} Let $\alpha \in  \overline{ME}(X)$ and $\sE$ and $\sF$ be 
$\alpha-$semi-stable torsion free sheaves on $X$.
Then $\sE \hat \otimes \sF$ is again $\alpha-$semi-stable. 
\end{theorem} 

First we treat a special case, which however is the crucial part of (5.1). 

\begin{proposition} Assume in the setup of (5.1) that $\alpha $ is geometric 
and $\sE$ and $\sF$ are locally free. Then 
$\sE \otimes \sF$ is $\alpha-$semi-stable. 
\end{proposition} 

\proof The proof is by adaption of the methods presented in [HL97] in the 
case that $\alpha = H^{n-1}$ with $H$ ample. \\
(1) The first step consists of the generalization of the Grauert-M\"ulich 
theorem 3.1.2 in [HL97] to the present situation. The conclusion in our case 
is the
same as if $\pi^*(\sE)$ were $(H_1, \ldots, H_{n-1})-$semi-stable. The proof 
is the same, but working with the 
complete intersection family $\tilde Z_s $ cut by the very ample divisors 
$H_i$ and their images $Z_s \subset X$. The Zariski open subset $X_0$ of [HL 
97] must be chosen
such that $X_0$ contains some $Z_s$. \\
(2) The next step carries over (3.1.6) of [HL97] to our situation. This is 
done verbatim, working on $\tilde X$. Here it is important to observe that 
$$ \mu_{\alpha}(\sE) = \mu_{(H_1,\ldots, H_{n-1})}(\pi^*(\sE)),$$ since 
$\sE$ is locally free. This would not work in case $\sE$ is only torsion 
free (or even reflexive),
since we cannot arrange $\sE$ to be locally free near the general $C.$ \\
(3) In the third step we generalize (3.2.8), Step 1, from [HL97]. This again 
works in the same way using the following easy fact. 
Suppose $Q$ is a torsion free sheaf
on $X$ of rank $r$. We define $\det Q = \bigwedge^r(Q)^{**}.$ If $\det(Q 
\vert C)$ is ample, then $(\det Q) \vert C$ is ample. \\
(4) The final step is the covering trick in Step 2 of (3.2.8). Here we 
choose $d_i \gg 0 $ and a branched covering
$$ \tilde f: \tilde X' \to \tilde X $$
such that $f^*(H_j) = \tilde H_j^{d_j}.$ 
Let $$ \tilde X' \buildrel {\pi'} \over {\to} X' \buildrel {f} \over {\to} X  
$$ 
be the Stein factorization of $ \pi \circ \tilde f.$ 
Let $\tilde C $ be the general curve cut out by the $H_i$ (supposed to be 
very ample) and, as always, $C = \pi(\tilde C)$ so that $\alpha = [C].$ 
The corresponding curves on $\tilde X'$ and $X'$ are denoted by $\tilde C'$ 
and $C'$. 
Then as in [HL97] we conclude that $f^*(\sE \otimes \sF)$ is 
$C'-$semi-stable, hence $f^*(C)-$semi-stable. Here we use the facts
$$ \pi'_*(\tilde f^*(\tilde C) = \pi'(\tilde f^{-1}(\tilde C)) = f^{-1}(C) = 
f^*(C) $$
and that $f^*(C)$ is a rational multiple of $C'.$ Now the semi-stability of 
$f^*(\sE \otimes \sF)$ w.r.t $f^*(C)$ implies the $C-$semi-stability of 
$\sE \otimes \sF$ (use the arguments of [HL97], (3.2.2)).

\qed

In order to deduce (5.1) from (5.2) we use the following lemma. 

\begin{lemma} Let $\sigma: \hat X \to X$ be a modification from the 
projective manifold $\hat X $. Let $\alpha \in  \overline{ME}(X)$ and
$\hat \alpha = \sigma^*(\alpha). $ Then 
\begin {enumerate} 
\item $\hat \alpha \in  \overline{ME}(\hat X)$
\item If $\sS$ is a torsion free sheaf on $X$ and $\hat S = 
\sigma^*(\sS)/{\rm tor},$ then $\mu_{\alpha}(\sS) = \mu_{\hat \alpha}(\hat 
\sS).$
\item If $\hat \sS$ is torsion free on $\hat X$ and $\sS = \sigma_*(\hat 
\sS),$ then  $\mu_{\alpha}(\sS) = \mu_{\hat \alpha}(\hat \sS).$
\item A torsion free sheaf $\sE$ on $X$ is $\alpha-$semi-stable if and only 
if $\sigma^*(\sE)/{\rm tor}$ is $\hat \alpha-$semi-stable.
\end{enumerate} 
\end{lemma}

\proof (1) We need to prove that $\hat D  \cdot \hat \alpha \geq 0$ for all 
pseudo-effective divisors $\hat D$ on $\hat X.$ Now 
$D = \sigma_*(\hat D)$ is again pseudo-effective. Hence
$$  \hat D  \cdot \hat \alpha \geq 0 = D \cdot \alpha \geq 0$$
proving (1). \\
(2) and (3) are simple calculations and (4) follows from (2) and (3). 
\qed

\proof {\it  of Theorem 5.1}. We proceed by induction on ${\rm rk}(\sE) + 
{\rm rk}(\sF).$ By Lemma 5.3 we may assume $\sE$ and $\sF$ locally free. 
If $\alpha$ is geometric, then Proposition 5.2 gives our claim. So assume 
that $\alpha$ is not geometric. Then we find a sequence $(\alpha_k)$
of geometric classes converging to $\alpha$ (in $N^1(X)$ or in $H^2(X)$); we 
can arrange $\alpha_t = \alpha + t(\omega_1 \wedge \ldots \wedge 
\omega^{n-1}$
with K\"ahler classes $\omega_i.$ \\
Suppose first that $\sE$ and $\sF$ are $\alpha-$stable, not just 
semi-stable. By Lemma 5.5 the bundles $\sE$ and $\sF$ are
$\alpha_k-$semi-stable for sufficiently large $k,$. Then by (5.2) $\sE 
\otimes \sF$ is $\alpha_k-$semi-stable for large $k$. Thus $\sE \otimes \sF$ 
is $\alpha-$semi-stable. \\
If $\sE$ and $\sF$ are $\alpha-$semi-stable but not both stable, we consider 
the saturated maximal destabilizing subsheaves $\sS \subset \sE$ and
$\sT \subset \sF$ with torsion free quotients $\sS'$ and $\sT'$ (possibly 
$\sS = \sE$ or $\sT = \sF$). 
By induction hypothesis $\sS \otimes \sT$ and $\sS' \otimes \sT'$ are 
semi-stable. Since
$$ \mu_{\alpha}(\sS) =  \mu_{\alpha}(\sS') =  \mu_{\alpha}(\sE) \eqno (1) $$
and
$$ \mu_{\alpha}(\sT) =  \mu_{\alpha}(\sT') =  \mu_{\alpha}(\sF) \eqno (2) $$
it follows easily that $\sE \otimes \sF$ is $\alpha-$semi-stable. Namely, 
tensor the exact sequence
$$ 0 \to \sT \to \sF \to \sT' \to 0 $$
by $\sS $ and $\sS'$ to deduce the semi-stability of $\sS \hat \otimes \sF$ 
and $\sS' \hat \otimes \sF$ and then tensor the the exact 
sequence 
$$ 0 \to \sS \to \sE \to \sS' \to 0$$
by $\sF$ to deduce the semi-stability of $\sE \hat \otimes \sF$. Here of 
course we need (1) and (2). 

\qed

\begin{corollary} 
Let $\alpha$ be a movable class on the projective manifold $X$ and $\sE$ and 
$\sF$ torsion free sheaves on $X$. Then 
$\mu^{\rm max}_{\alpha}(\sE \hat \otimes \sF) =  \mu^{\rm max}_{\alpha}(\sE) 
+ \mu^{\rm max}_{\alpha}(\sF). $  

\end{corollary}

\proof Let $\sS \subset \sE $ and $\sT \subset 
\sF$ be the maximal destabilizing sheaves. Since $\sS \hat \otimes \sT$ is 
$\alpha-$semi-stable by Theorem 5.1, we obtain
$$\mu^{\rm max}_{\alpha}(\sE \hat \otimes \sF) \geq \mu_{\alpha}(\sS \hat 
\otimes \sT).$$
Since $\mu^{\rm max}_{\alpha}(\sE) = \mu_{\alpha}(\sS),$ and analogously for 
$\sF$ and $\sT,$ we conclude for one inequality.
To establish the other, we must show that $\sS \hat \otimes \sT$ is maximal 
destablisizing for $\sE \hat \otimes \sF.$ 
This is an easy exercise using the exact sequences already used in the proof 
of (5.1) and the HN-filtration.

\qed

\begin{lemma} Let $X$ be a projective manifold of dimension $n$ and $\sE$ a 
reflexive sheaf over $X$. Let $ \alpha \in {\overline {ME}}(X).$  Choose 
K\"ahler classes $\omega_i$ and set
$$ \alpha_t = \alpha + t(\omega_1 \wedge \ldots \wedge \omega^{n-1})$$
for $t \in \bR_+.$
Assume that $\sE$ is $\alpha-$stable. Then there exists a sequence $(t_j) $ 
converging to $0$ such that $\sE$ is $\alpha_{t_j}-$semi-stable. 
\end{lemma}

\proof We assume to the contrary that $\sE$ is $\alpha_t-$unstable for $0 < 
t \leq t_0.$ 
\vskip .1cm \noindent
(1) Let first $(s_i)$ be a sequence in $(0,t_0] $ converging to $s_0 > 0$ 
such that some fixed saturated subsheaf $\sS \subset \sE$ is maximal 
destabilizing for all
$\alpha_{s_i}$. Let $\sS_0$ be the maximal $\alpha_{s_0}-$destabilizing 
subsheaf. Then either $\sS = \sS_0$ or $\sS \subset \sS_0$ with 
${\rm rk} \sS < {\rm rk} \sS_0.$ \\
In fact, in the obvious notation we have 
$$ \mu_{s_i}(\sS) \geq \mu_{s_i}(\sS_0) $$
for all $s_i.$ Using the assumption that $\sS$ is $\alpha_{s_0}-$maximal 
destabilizing, we obtain in the limit
$$ \mu_{s_0}(\sS) = \mu_{s_0}(\sS_0). $$
Therefore $\sS \subset \sS_0.$ Suppose that both (reflexive) sheaves have 
the same rank, but are not equal. Then 
$\det S_0 = \det S + D$ with $D$ effective non-zero, which clearly 
contradicts the last equality, since $\alpha_{s_0}$ is inside the interior 
of the
movable cone. 
\vskip .1cm \noindent  (2) From (1) we deduce that after possibly shrinking 
$t_0,$ the set of positive real numbers $t \leq t_0$ 
for which a given subsheaf $\sE \subset \sE$ is $\alpha_t-$maximal 
destabilizing, is closed. 
\vskip .1cm \noindent
(3) Now let $\sS_t$ be the maximal destabilizing subsheaf for $\alpha_t.$ We 
will be finished if we can prove that $\sS_t$ is independent of $t$ at 
least for small positive $t.$ We start with $\sS_0 = \sS_{t_0}$ and choose, 
using (2), $t_1$ minimal such that $\sS_0$ is still maximal destabilizing
for $\alpha_{t_1}.$ 
Consider $t < t_1$ near to $t_1$ and consider the induced map
$$ \phi_{0,t}: \sS_0 \to \sE/\sS_t.$$
Our aim is to show that 
$$\phi_{0,t} = 0 \eqno (*)$$ 
for $t$ sufficiently to $t_1$. 
Given (*), we conclude that $\sS_0 \subset  \sS_t$ and as in (1), $\sS_t$ 
will have larger rank. Proceeding then with some $t < t_1$ near to $t_1$, 
this jumping behaviour can happen only finitely many times, so that after 
finitely many steps, we are able to conclude that $\sS_t$ is independent
of $t.$ \\
It remains to prove (*). First observe that since
$$ \mu_{t_1}(\sS_0) > \mu_{t_1}(\sE),$$
the same will be true for $t$ near to $t_1.$ 
Since $\mu_{t}(\sS_t) > \mu_{t}(\sE)$ by assumption, we obtain
$$ \mu_t(\sS_0) > \mu_t(\sE/\sS_t). \eqno (**) $$
(3a) If $\sS_0$ is $\alpha_{t_1}$-stable, then by induction on the rank we 
may assume that $\sS_0$ is $\alpha_t$-stable for $t$ near $t_1$, hence
from (**), the vanishing (*) follows. \\
(3b) So suppose that $\sS_0$ is $\alpha_{t_1}-$semi-stable but not stable. 
Let $\sS'_0 \subset \sS_0$ be
the maximal $\alpha_{t_1}-$destabilizing subsheaf. Then $\sS_0/\sS'_0$ is 
again $\alpha_{t_1}-$semi-stable. 
The restricted map $\phi_{0,t} \vert  \sS'_0$ vanishes by the same arguments 
as in (3a). Thus we obtain a map
$$ \phi'_{0,t}: \sS_0/\sS'_0 \to \sE/\sS_t.$$ 
Now we go back to the beginning of (3b) substituting $\sS_0$ by 
$\sS_0/\sS'_0$ and $\sE$ by $\sE/\sS_t.$ in other words, we proceed 
inductively
and conclude that $ \phi'_{0,t} = 0,$ hence (*) holds.
\qed
 
\newpage

\section{Appendix: an alternative proof of theorem 5.1.}



  We give here an alternative proof of Theorem 5.1.. 
  This proof needs no adaptation of the Grauert-M\"ullich Theorem. 
  The main ingredient will be the Kobayashi-Hitchin correspondence for non-K\"ahler polarizations 
  which was established by J.Li and S.-T. Yau.


We shall use the notations and the definitions of the main paper. In particular $X$ will be a complex projective manifold of dimension $n\ge 2$. Following 
\cite{DPS96} we shall denote by $N_{amp}$ the interior of the closed cone $\overline{ME(X)}$ generated by movable curves, see also  \cite{BDPP04}. 
It is easy to check using \cite{BDPP04} that geometric classes of curves belong to $N_{amp}$.

The following proposition replaces Proposition 5.2.

\begin{proposition} Let $\alpha$ be a class in $N_{amp}$ and $\mathcal E$ and $\mathcal F$ two $\alpha$-polystable locally free sheaves. Then 
$\mathcal E\otimes \mathcal F$ is again 
 $\alpha$-polystable.
\end{proposition} 
\begin{proof}
   We start by a Hahn-Banach argument 
   and show the existence of  a smooth positive definite form $u$ of bidegree $(n-1,n-1)$ with
    $\partial\overline{\partial} u=0$ and such that 
    the slope of a holomorphic vector bundle with respect to $\alpha$ is computed by 
    $$\mu_{\alpha}(E)=\frac{\int c_1(E,h) \wedge u}{\rank E},$$
    where $c_1(E,h)$ is the first Chern form of $E$ computed with respect to some hermitian metric $h$ in the fibers of $E$,
    cf. \cite{DPS96} Theorem 4.1.
    
    Let indeed 
    $\mathcal D^+_{1,1}$ be the cone of positive currents inside the space 
    $\mathcal D'_{1,1} $ of currents of bidegree $(1,1)$.
    For any choice of a  positive definite smooth  $(1,1)$-form $\eta$ the set 
    $\mathcal D^+_{(1,1),\eta}=\{ T\in \mathcal D^+_{1,1} \ | \ \int_X T\wedge\eta^{n-1}=1\}$ 
    is compact for the weak topology on
    $\mathcal D'_{1,1}$, see \cite{D} III.1.23. The vector subspace 
    $V=\{ T\in \mathcal D'_{1,1} \ | \ \d T=0, \ [T]\alpha=0\}$
    is closed and disjoint from $\mathcal D^+_{(1,1),\eta}$ by \cite{BDPP04}. 
    (Notice that $\alpha$ belongs also to the interior
    of the cone of movable classes $\mathcal{M}\subset H^{n-1,n-1}_{\bR}(X)$.)
    Thus there exists a continuous linear functional which is positive on $\mathcal D^+_{(1,1),\eta}$ and vanishes on $V$.
    This is given by a smooth positive definite form $u$ of bidegree $(n-1,n-1)$ which is also
    $\partial\overline{\partial}$-closed since $\partial\overline{\partial}\mathcal D'_{0,0}\subset V$. 
    Moreover a renormalization
    of $u$ by a positive factor makes 
    $\alpha$ and $u$ to be equal as linear functionals on 
    $H^{1,1}(X)_{\bR}$ since they have the same kernel and are both positive on K\"ahler classes.

    Next we obtain a  $(n-1)$-st root $\omega$ of $u$ in the following way. First notice that 
     $$(i\sum_{1\le i,j\le n}a_{ij}\d z_i\wedge\d \overline z_j)^{n-1}=
     (n-1)!i^{(n-1)^2}\sum_{1\le i,j\le n}(-1)^{i+j}c_{ji}\hat{\d z_i}\wedge\hat{\d \overline z_j},\leqno{(1)}$$
     where 
     $c_{ij}$ denotes the cofactor of $a_{ij}$ in the matrix $A=(a_{ij})_{1\le i,j\le n}$,
     $\hat{\d z_i}=\d z_1\wedge...\wedge\d z_{i-1}\wedge\d z_{i+1}\wedge...\wedge\d z_n$ and 
     $\hat{\d \overline z_j}
     =\d \overline z_1\wedge...\wedge\d\overline z_{j-1}\wedge\d\overline z_{j+1}\wedge...\wedge\d \overline z_n$.
     The relation $\trans{C}A=\det (A)I_n$ for the cofactor matrix $C=(c_{ij})_{1\le i,j\le n}$ implies 
     $$A=\sqrt[n-1]{\det(C)}\trans{C}^{-1}$$
     in case $A$ is positive definite. Moreover, given a positive definite
     matrix $C$, one obtains a unique positive definite solution $A$ of the equation (1). 
     
     Then $\omega$ is the $(1,1)$-form associated to a Gauduchon metric on $X$
     and  
     $$\mu_{\alpha}(E)=\mu_{\omega}(E)=\frac{\int  c_1(E,h)\wedge \omega^{n-1}}{\rank E},$$
     for $E$ and $h$ as before.
     By \cite{LY87} the Kobayashi-Hitchin correspondence holds in this case, 
     thus the polystability of a holomorphic vector bundle $E$ with respect to $\omega$ 
     is equivalent to the existence of a Hermite-Einstein metric with respect to the polarization $\omega$ again.
     But the tensor product of Hermite-Einstein vector bundles is also Hermite-Einstein and the Proposition is proved.
\end{proof}

We can now prove:

\begin{theorem}
Let $\alpha\in\overline{ME(X)}$ and $\mathcal E$ and $\mathcal F$ two $\alpha$-semi-stable torsion free sheaves. Then 
$\mathcal E\hat{\otimes}\mathcal F$ is $\alpha$-semi-stable.
\end{theorem}
\begin{proof}
The proof follows basically the strategy of the main paper and uses the above proposition instead of 
Proposition 5.2.

We will make induction on $\rank\mathcal E+\rank\mathcal F$.

By  Lemma 5.3 we may assume that $\mathcal E$ and $\mathcal F$ are locally free.
 
We may also assume that $\mathcal E$ and $\mathcal F$ are $\alpha$-stable: if $\mathcal E$ is $\alpha$-semi-stable but not $\alpha$-stable there is some coherent subsheaf of lower rank $\mathcal E_1\subset\mathcal E$ with torsion free quotient such that 
$\mu_{\alpha}(\mathcal E_1)=\mu_{\alpha}(\mathcal E)=\mu_{\alpha}(\mathcal E/\mathcal E_1)$. By the induction hypothesis we get $\mathcal E_1\otimes \mathcal F$ and 
$(\mathcal E/\mathcal E_1)\otimes \mathcal F\cong (\mathcal E\otimes \mathcal F)/(\mathcal E_1\otimes \mathcal F)$ $\alpha$-semi-stable of the same slope hence also 
$\mathcal E\otimes \mathcal F$ will be  $\alpha$-semi-stable.

By Lemma 5.5 we may approximate the class $\alpha$ by a sequence $(\alpha_k)_{k\in \bN}$  of classes in $N_{amp}$
such that both $\mathcal E$ and $\mathcal F$ are  $\alpha_k$-semi-stable for all $k\in\bN$.  Our proposition combined with a use of the induction
hypothesis as above implies that  $\mathcal E\otimes \mathcal F$ is $\alpha_k$-semi-stable for all $k\in\bN$.
Thus   $\mathcal E\otimes \mathcal F$ will be $\alpha$-semi-stable as well.

\end{proof}

\vspace{1cm}
\small
\begin{tabular}{lcl}
Fr\'ed\'eric Campana && Thomas Peternell \\D\'epartement de
Math\'ematiques && Mathematisches Institut \\
Universit\'e de Nancy && Universit\" at Bayreuth\\
F-54506 Vandoeuvre-les-Nancy, France && D-95440 Bayreuth, Germany \\
frederic.campana@iecn.u-nancy.fr && thomas.peternell@uni-bayreuth.de\\
\end{tabular}

\end{document}